\documentclass[12pt]{article}
\usepackage[cp1251]{inputenc}
\usepackage{epic}
\usepackage{eepic}
\usepackage{graphics}
\usepackage{amsfonts}
  \usepackage{amsthm}
 \usepackage{amsmath}
 \usepackage{amscd}
 \usepackage{amssymb}
 \usepackage{latexsym}
\usepackage{graphicx}

 \numberwithin{equation}{section}
 \newtheorem{thm}{Theorem}
 
 \newtheorem{prop}{Proposition}
 \theoremstyle{definition}

\title{ Non-integrability of a system with  the Dyson Potential}

\author{Georgi Georgiev }
\date{}

\begin{document}

\maketitle

\begin{abstract}
In this paper it is shown that the Hamiltonian system with Dyson potential is analytical non-integrable and formal non-integrable.  The approach is based on the following: when a system has a family of periodic solutions around an equilibrium and if the period function is infinitely branched then the system has non additional analytic first integral. We prove formal non-integrability  using Ziglin-Moralez-Ruiz-Ramis theory.
\end{abstract}



We study the Hamiltonian system of $n$ interacting particles of equal mass with a Hamiltonian 
\begin{equation}
\label{1.1}
H=\frac{1}{2} \sum_{i=1}^{n} y_i^2+\sum_{i<j} V(x_i-x_j),\, V(x)=-\log|\sin x|.
\end{equation}
Here $x_1,\, x_2,\, \dots x_n$ are the coordinates of particles, $y_1,\, y_2, \, \dots y_n$ are their momenta, and $V$ is the potential energy in Dyson type.  
The system with this potential was studied in \cite{Dys}, where the statistical properties of the energy levels of  one dimensional Coulomb's gas are investigated.There is a connection between system of point vortices and regarding system in \cite{CalPer} by Calogero and Perelomov. The equilibria position of the system (\ref{1.1}), $x_k^0=x_0+\frac{\pi k}{n},\, k=1,\, 2, \, \dots n, \, x_0\in {\bf R}$ determines stationary collinear choreography on the sphere-point vortices are located in the equatorial plane that uniformly rotates around an axis lying in this plane.

The function $V$ is $2\pi$-periodic (and even $\pi$-periodic), so we can assume that the particles move in circles. The Hamiltonian system always has two integrals
\begin{equation*}
H,\, F=\sum_{i=1}^{n} y_i.
\end{equation*}
The question of integrability of  (\ref{1.1}) is studied by Calogero and Perelomov \cite{CalPer}. If $V$ is a non-constant analytic periodic function without singularities, then the system (\ref{1.1}) can not be integrable for $n\ge 3$ ( Kozlov \cite{Koz1}). Unfortunately Dyson potential  has a real logarithmic singularity.  Borisov and Kozlov in \cite{BorKoz}  had proved  that the system in case $n=3$ is non-integrable in analytical first integrals.
We consider the same case  $ n=3$, but  the approach is different, this is the first non-trivial case and we prove that the system admits only $F$ as a holomorphic first integral in a complex domain.

 Let us make in (\ref{1.1})  (in our  case $n=3$) canonical change of variables 
\begin{eqnarray*}
y_1 & = &p_1+p_3,\, y_2=-p_1+p_2+p_3, \, y_3=-p_2+p_3,\\
q_1 & =&  x_1-x_2,\, q_2=x_2-x_3, \, q_3=x_1+x_2+x_3,
\end{eqnarray*}
and we obtain
\begin{equation}
\label{1.5}
H=p_1^2-p_1p_2+p_2^2+\frac{3}{2}p_3^2-\log|\sin {q_1}|-\log|\sin {q_2}|-\log|\sin {(q_1+q_2)}|.
\end{equation}
The system has another integral
\begin{equation}
\label{1.51}
F=3p_3
\end{equation}
and has a stable equilibrium  $p_1=p_2=0$, $q_1=q_2=\displaystyle{\frac{\pi}{3}}$. 


The variable $q_3$ is cyclic and it  reduces the system (with $p_3=const$) to a system with 2-degrees of freedom  with Hamiltonian 
\begin{equation}
\label{1.66}
H_{reg}=p_1^2-p_1p_2+p_2^2-\log|\sin {q_1}|-\log|\sin {q_2}|-\log|\sin {(q_1+q_2)}|.
\end{equation}

We need  some basic definitions about Hamiltonian systems.

Let $H$ is a smooth real-valued function of $2n$ real variables $(p,\, q),\, p,\,q\, \in {\bf R^n}$. Let us also  assume that $dH(0)=0$ where $0$ is an equilibrium point for the Hamiltonian system  $X_H$ (with $n$ degrees of freedom), given by

\begin{equation*}
\frac{dq}{dt}=\frac{\partial H}{\partial p},\,\, \frac{dp}{dt}=-\frac{\partial H}{\partial q}.
\end{equation*}
We often write the Hamiltonian systems in the form

\begin{equation*}
\dot{x}=X_{H}(x),\, x\in {\bf R^{2n}},
\end{equation*}
where $X_{H}$ is the flow.
The system is called Liouville - Integrable near $0$ if there exists $n$ functions in involution $f_1=H,\, f_2,\, \dots f_n$, defined around $0$, are functionally independent. The Poisson bracket of $f$ and $g$ are
 \begin{equation*}
\{f,\,g\}=X_{f} (g) =\sum\frac{\partial f}{\partial p}\frac{\partial g}{\partial q}-\frac{\partial f}{\partial q}\frac{\partial g}{\partial p}=-\{g,\, f\}.
\end{equation*}
We call that  functions $f$ and $g$ are in involution if the Poisson bracket is commuteve.
This means that $df_1,\, df_2,\, \dots df_n$ are linearly independent around the equilibrium $0$ and $f_j=const$ for all $j$ define smooth submanifolds, these level manifolds are invariant under $X_{f_i}$. We have $X_{f_j}f_{j}=0$ and $[X_{f_j},\,X_{f_k}]=X_{\{f_{j},f_{k}\}}=0$ - these flows commute. The compact and connected component of $M_c:=\{ f_{j}=c_j,\, j=1,\, \dots , n\}$ is diffeomorphed to a torus.

We call that the system with Hamiltonian $H$ formally integrable if there exist formal power series $\tilde f_1,\, \tilde f_2,\ \dots \tilde f_n$ in involution, where  $\tilde f_1,\, \tilde f_2,\ \dots \tilde f_n$ are functionally independent and Taylor expansion $\tilde H$ of $H$ is a formal power series in $\tilde f_1,\, \tilde f_2,\ \dots \tilde f_n$. An asymptotic behavior near equilibrium is like an integrable system. The formal integrability gives information about the flow. The functional independence of $\tilde f_1,\, \tilde f_2,\ \dots \tilde f_n$ is  stronger than the independence of smooth functions $f_1,\, f_2,\, \dots f_n$ of which $\tilde f_1,\, \tilde f_2,\ \dots \tilde f_n$ are the Taylor series, because formal independence leads to a functional independence of a finite part of Taylor expansions.

Here we study the formal non-integrability of the system (\ref{1.1})  in case $n=3$.

Our first  aim is the following

\begin{thm}
\label{th1}
a) The  system with Hamiltonian (\ref{1.66}) 
is not  integrable by means of analytical first integral;

b)  The system with Hamiltonian(\ref{1.66}) 
is not  formal integrable.
\end{thm}

The motivation for proving formal non-integrability I received from the remarcable paper of J. J. Duistermaat \cite{Duis}.  The Hamiltonian system with Dyson potential has similar structure. The proof of a) is different from \cite{BorKoz}.

{\bf Proof a):}
 The proof  of  this  is based on three propositions. The first - shows that there  exists family of periodical solutions near equilibrium. The second proposition investigates the behavior on period function of these periodical solutions around the equilibrium. The third - proves that if we assume an existence of an additional first integral- it is a constant.

\begin{prop}
\label{pr1}
On the manifold $P:=\{p_1=p_2,\, q_1=q_2\}$, invariant under $X_H$ the system with Hamiltonian (\ref{1.66}) exists a family of periodical solutions around equilibrium $p_1=p_2=0$, $q_1=q_2=\displaystyle{\frac{\pi}{3}}$.
\end{prop}

{\bf Proof:}
Let we take  $q(t):=q_1(t)=q_2(t),\, p(t):=p_1(t)=p_2(t)$ and $\displaystyle{\frac {d\, q(t)}{dt}}=p(t)$ in our system, then

\begin{equation*}
\dot p   =  \cot q +\cot 2q , 
\end{equation*}
and we have
\begin{equation*}
(\dot q )^2=\log(\sin 2q)+2\log(\sin q)+E.
\end{equation*}

This is a conservative system with  convex potential 
\begin{equation}
\label{1.10}
\tilde V(q)=-\log(\sin 2q)-2\log(\sin q),
\end{equation}
and we know that for these predictions the system has a  periodical solution around  the equilibrium for fixed energy $H_{reg}=E$ (see Figure 1). 

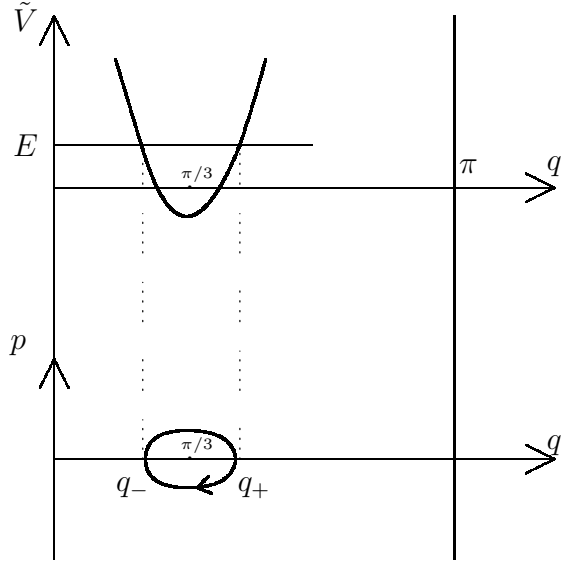
\begin{figure}[h]
\begin{center}
 \setlength{\unitlength}{1.9mm}
 \begin{picture}(35,45)(0,0)
 \linethickness{0.4pt}

       \qbezier(0,34)(15,34)(18,34)  
       \qbezier(0,31)(15,31)(35,31)  

       \put(-2,34){\makebox(0,0)[cc]{$E$}}
       \put(35,32.5){\makebox(0,0)[cc]{$q$}}

   \put(29,32.5){\makebox(0,0)[cc]{$\pi$}}

       \qbezier(35,31)(34,30.5)(33,30)   
       \qbezier(35,31)(34,31.5)(33,32)   
       \qbezier(0,10)(0,38)(0,43)    
       \qbezier(0,43)(0.5,42)(1,41)  
       \qbezier(0,43)(-0.5,42)(-1,41)    
       \put(-2,43){\makebox(0,0)[cc]{$\tilde{V}$}}


 \qbezier(28,5)(28,38)(28,43)

\dashline{3}[0.7](6.2,12)(6.2,34)

\dashline{3}[0.7](13,12)(13,34)

  \put(5.5,10){\makebox(0,0)[cc]{$q_{-}$}}
  \put(14,10){\makebox(0,0)[cc]{$q_{+}$}}

       \qbezier(0,12)(15,12)(35,12)  
       \put(35,13){\makebox(0,0)[cc]{$q$}}


       \qbezier(35,12)(34,12.5)(33,13)   
       \qbezier(35,12)(34,11.5)(33,11)   
       \qbezier(0,5)(0,12)(0,19)    
       \put(-2.5,20){\makebox(0,0)[cc]{$p$}}
       \qbezier(1,17)(0.5,18)(0,19)  
       \qbezier(-1,17)(-0.5,18)(0,19)    

 \linethickness{0.8pt}
       \qbezier(6.1,34)(9,24)(13,34)
       \qbezier(13,34)(14,37)(14.8,40)
       \qbezier(4.3,40)(5.2,37)(6.1,34)

 \linethickness{0.8pt}

    \qbezier(6.4,12)(6.4,14)(9.3,14)
    \qbezier(9.3,14)(12.5,14)(12.7,12)
    \qbezier(9.3,10)(12.5,10)(12.7,12)
    \qbezier(10,10)(10.4,10.4)(10.8,10.8)
    \qbezier(10,10)(10.5,9.8)(11,9.6)
    \qbezier(6.4,12)(6.4,10)(9.3,10)

\put(9.5,31){\makebox(0,0)[cc] {${\cdot}$}}
\put(10,32){\makebox(0,0)[cc] {${\scriptscriptstyle \pi/3}$}}


\put(9.5,12){\makebox(0,0)[cc] {${\cdot}$}}
\put(10,13){\makebox(0,0)[cc] {${\scriptscriptstyle \pi/3}$}}

 \end{picture}
$\
$
\end{center}
\caption{\hspace{-0.5pt} Periodical solution.}
\end{figure}

\begin{prop}
\label{pr2}
The period function has expression $T(c)=\log{\eta(c) }+\Phi(c)$, where $c=\displaystyle{\frac{1}{2e^E}}$, $\eta(c)=\epsilon (B(c))\delta (B(c))$ with
\begin{equation}
\label{B}
B(c)=\frac{16.(\frac{2}{3})^{1/3}.c^2}{(9c^2-\sqrt{3}\sqrt{27c^4-256c^6})^{1/3}}+2(\frac{2}{3})^{1/3}.(9c^2-\sqrt{3}\sqrt{27c^4-256c^6})^{1/3}.
\end{equation}

$\Phi(c)$ is an analytical function  of the variable $c$ and 
 \begin{equation}
 \label{p_1}
 r_{1,2}=\frac{1}{2}-\frac{1}{2}\sqrt{1+B}\pm \frac{1}{2}\sqrt{2-B+\frac{2}{\sqrt{1+B}}},
 \end{equation}
 $\epsilon=\displaystyle{\frac{1}{2}(\arccos{r_1}-\frac{\pi}{3})}$, $\delta=\displaystyle{\frac{1}{2}(\frac{\pi}{3}-\arccos{r_2})}$.
\end{prop}
{\bf Proof:}
The period in this solution  is
\begin{equation}
\label{1.7}
T=2\int_{q_-}^{q_+}\frac{dq}{\sqrt{E-\tilde V(q)}},
\end{equation}
where $q_{-}$ and $q_+$ are the roots of $E-\tilde V(q)=0$, where $\tilde V(q)$ is (\ref{1.10}).

For equilibrium point $q=\displaystyle{\frac{\pi}{3}}$ we have $\tilde{V}(\displaystyle{\frac{\pi}{3}})=-3\log(\displaystyle{\frac{\sqrt{3}}{2}})$ and let we fix $q_{+}=\displaystyle{\frac{\pi}{3}+\epsilon} $ and $q_{-}=\displaystyle{\frac{\pi}{3}-\delta} $, for $\epsilon >0$ and $\delta >0$. We have
\begin{equation*}
A=\frac{2}{\sqrt{E-\tilde V(q)}}=\frac{2}{\sqrt{-3\log(\displaystyle{\frac{\sqrt{3}}{2}})+\log(\sin 2q)+2\log(\sin q)}}.
\end{equation*}
The expansion of $A$ near $\displaystyle{\frac{\pi}{3}}$ is
\begin{equation*}
 A=\frac{1}{|q-\frac{\pi}{3}|}-\frac{1}{3\sqrt{3}}-\frac{2}{9}|q-\frac{\pi}{3}|+O(|q-\frac{\pi}{3}| ^2). 
\end{equation*}
For the period we obtain 
\begin{eqnarray*}
T & = & \int_{\frac{\pi}{3}-\delta}^{\frac{\pi}{3}+\epsilon}\frac{dq}{\sqrt{E-\tilde{V}(q)}}
   =  \int_{\frac{\pi}{3}-\delta}^{\frac{\pi}{3}+\epsilon}\left(\frac{1}{|q-\frac{\pi}{3}|}-\frac{1}{3\sqrt{3}}-\frac{2}{9}|q-\frac{\pi}{3}|+\dots \right)dq\\
   & = & \log(\epsilon)+\log(\delta)+\Phi(\epsilon ,\, \delta)=\log(\epsilon .\delta)+\Phi(\epsilon ,\, \delta)=\log{\eta}+\Phi(\epsilon,\, \delta),
\end{eqnarray*}
where $\eta=\epsilon\delta$ and $\Phi(\epsilon ,\, \delta)$ is an analytical function on each of the variables $\epsilon$ and $\delta$. Next we find $\epsilon$ and $\delta$. We have 
$$E=\tilde{V}(q)=-\log(\sin 2q)-2\log(\sin q)=-\log(2.(\sin q)^3.\cos q),$$
and if we get $q_{+}=\displaystyle{\frac{\pi}{3}+\epsilon}$ ( or $q_{-}=\displaystyle{\frac{\pi}{3}-\delta}$), then - $\displaystyle{(\sin(\frac{\pi}{3}+\epsilon))^3.\cos(\frac{\pi}{3}+\epsilon)=\frac{1}{2e^E}=c}$. If we put $r=\cos(\displaystyle{\frac{2\pi}{3}}+2\epsilon)$, we obtain equation $(1-r)^2(1-r^2)=16c^2$. The real roots of this equation are (\ref{p_1})  with (\ref{B}).  
We find $\epsilon=\displaystyle{\frac{1}{2}(\arccos{r_1}-\frac{\pi}{3})}$, $\delta=\displaystyle{\frac{1}{2}(\frac{\pi}{3}-\arccos{r_2})}$.

\begin{prop}
\label{pr3}
The system with Hamiltonian (\ref{1.66}) does not possess any  additional holomorphic first integral. 
\end{prop}
{\bf Proof:}
 It is important for the proof that there is a family of solutions on $P$, of the Hamiltonian system (\ref{1.66}) on the hypersurface $H=E$ are periodic. If $T$ is the period function then complex continuation of the manifolds $T=const$ turns out to be infinitely branched. This excludes the existence of a nontrivial analytic integral on any open subset of the complex domain where this infinite branching is true. 
Further, we need to show that if $G$ is a smooth function on open set $U$ such that $V=U\cap (H=E)$ is $X_H$ invariant, $\{H,\, G\}=0$ and derivative $d\{H,\, G\}=0$ on $V$, then $G$ is a function of $H$ and $T$ on $V$.
 The function $T(x)$, $x\in\tilde{U}\subset U$ is real analytical and it is not constant  $P$, it means  that $dT\ne 0$ an open dense subset $\Delta$  of the manifold $P$, $G$ is invariant under flows $X_H$ and $X_T$. We have $X_G T=0$ on $V$. Further we regard as Hamiltonian flow on $P$ of any smooth extension $\tilde T$ on $T$ on open neighborhood of $P$, the flow is independent of the choice of the extension. On the open set $\Delta\cap V$, $dT$ and $dH$ are linearly independent that is why the flows  $X_H$ and $X_T$ walk around one-dimensional submanifold $P$, $H=E$ and $T=const$. Therefore on $\Delta\cap V$, $G$ is locally constant on  $H=E$, and $T=const$. That is $G$ function of $H$ and $T$ on the connected component of $\Delta\cap U\cap P$.
 
 Let suppose that $G$ is analytic on $U$ and it has complex analytic extension on $\tilde U$ of $U$. If $G$ is not functionally depenended on $H$ on $P$ then the  manifold  $P$  and $G=const$ extend to closed complex analytic manifold on $\tilde U$, which  coincides with complex analytic continuations of  $H=E$  and $T=const$. If the analytic continuations have infinite branching near $x\in \tilde U$, then we have a contradiction, so $G$ has to be a function of $H$ on  $P$. Now we use $H$ as a coordinate near x, we can write $G=G_0+G_1H$, $G_0$ is a function of $H$ and $G_1$ is analytic. $G_1$ commutes with $H$ that is why $G_1$, this gives us that $G_1$ is a function of $H$ on $P$.  We obtain that $G$ is a function of  $H$ near $x$. By the analytic continuation this is true in the connected component of $x\in U$. 
 
{\bf Proof b):}
Let us go back to the formal non-integrability. First we will regard the case $K=H_2 +H_3$ it is Taylor expansion to a degree 3 with change of variables $\tilde{p_1}=p_1$, $\tilde{p_2}=p_2$, $\tilde{q_1}=q_1-\frac{\pi}{3}$ and $\tilde{q_2}=q_2-\frac{\pi}{3}$ (we move equilibrium to $0$)
\begin{equation}
K=\tilde{p_1}^2-\tilde{p_1}\tilde{p_2}+\tilde{p_2}^2+\frac{4}{3}\tilde{q_1}^2+\frac{4}{3}\tilde{q_1}\tilde{q_2}+\frac{4}{3}\tilde{q_2}^2+
\frac{4}{9}\sqrt{3}\tilde{q_1}^2\tilde{q_2}+\frac{4}{9}\sqrt{3}\tilde{q_1}\tilde{q_2}^2.
\end{equation} 
Let us remove tildes and we obtain the Hamiltonian system
\begin{eqnarray*}
\dot{q_1} & = & 2p_1-p_2\\
\dot{q_2} & = & -p_1+2p_2\\
\dot{p_1} & = & -\frac{8}{3}q_1-\frac{4}{3}q_2-\frac{8}{9}\sqrt{3}q_1q_2-\frac{4}{9}\sqrt{3}{q_2}^2\\
\dot{p_2} & = & -\frac{4}{3}q_1-\frac{8}{3}q_2-\frac{8}{9}\sqrt{3}q_1q_2-\frac{4}{9}\sqrt{3}{q_1}^2
\end{eqnarray*}


We use  the Theory of Morales-Ruiz- Ramis (see \cite{MR1}  for details) reducing  the system to a Normal Variations Equations (NVE) near a non-trivial partial solution. We find  a partial solution for $p=p_1=p_2$, $q=q_1=q_2$,
and we have
\begin{equation*}
\dot{q}^2=-\frac{8}{9}\sqrt{3}q^3-4q^2+h,
\end{equation*}
with solution $\phi(t)=\displaystyle{-\frac{\sqrt{3}}{2}-\frac{3\sqrt{3}}{2}\wp(t,g_2,g_3)}$, here $\wp$ is Weierstrass p-function with  $g_2=\displaystyle{\frac{4}{3}}$,  $g_3=-\displaystyle{\frac{4(h-2)}{27}}$.   
Let we put $\eta_1=dp_1$, $\eta_2=dp_2$, $\xi_1=dq_1$ and $\xi_2=dq_2$ then
\begin{eqnarray*}
\dot{\xi_1} & = & 2\eta_1-\eta_2\\
\dot{\xi_2} & = & -\eta_1+2\eta_2\\
\dot{\eta_1} & = & -(\frac{8}{3}+\frac{8}{9}\sqrt{3}. \phi )\xi_1 - (\frac{4}{3}+\frac{16}{9}\sqrt{3}. \phi )\xi_2\\
\dot{\eta_2} & = & -(\frac{4}{3}+\frac{16}{9}\sqrt{3}. \phi )\xi_1 - (\frac{8}{3}+\frac{8}{9}\sqrt{3}. \phi )\xi_2.
\end{eqnarray*}
If we get $\xi=\xi_1-\xi_2$, then we find an equation for NVE
$\ddot{\xi}=(-\frac{8}{3}+4\wp(t,g_2,g_3))\xi $.
This is a Lame-equation with $A=4$, $B=-\frac{8}{3}$ and in this case we have non Lame-Hermite solutions $A=n(n+1)\ne 4$ for $n\in {\bf Z}$.  We have non  Briochi- Halphen- Crowford solutions $n+1/2\notin {\bf N}$  and we have non the Baltassarri solutions $n+1/2\notin \frac{1}{3}{\bf Z}\cap\frac{1}{4}{\bf Z}\cap\frac{1}{5}{\bf Z}-{\bf Z}$  -  the identity component of its Galois group is non-commutative (see \cite{ChrGeor} for details). 

The theory says that if the identity component of differential Galois  group is non-commutative, then the system is not meromorphic  integrable (see \cite{MR1}). This proves that in the case $H_2+H_3$ there is non additional meromorphic (holomorphic) first integral.

Let us consider the case $H_2+H_3+H_4$: the Hamiltonian  is
\begin{eqnarray}
\label{L}
\Lambda & = & \tilde{p_1}^2-\tilde{p_1}\tilde{p_2}+\tilde{p_2}^2+\frac{4}{3}\tilde{q_1}^2+\frac{4}{3}\tilde{q_1}\tilde{q_2}+\frac{4}{3}\tilde{q_2}^2+
\frac{4}{9}\sqrt{3}\tilde{q_1}^2\tilde{q_2}\nonumber \\
& + & \frac{4}{9}\sqrt{3}\tilde{q_1}\tilde{q_2}^2+\frac{4}{9}\tilde{q_1}^4+\frac{4}{9}\tilde{q_2}^4
+\frac{8}{9}\tilde{q_1}^3 \tilde{q_2}+\frac{8}{9}\tilde{q_1} \tilde{q_2}^3+\frac{4}{3}\tilde{q_1}^2 \tilde{q_2}^2.
\end{eqnarray}
 We ignore the tildes and we get the system

\begin{eqnarray*}
\dot{q_1} & = & 2p_1-p_2\\
\dot{q_2} & = & -p_1+2p_2\\
\dot{p_1} & = & -\frac{8}{3}q_1-\frac{4}{3}q_2-\frac{8}{9}\sqrt{3}q_1q_2-\frac{4}{9}\sqrt{3}{q_2}^2\\
          & - & \frac{16}{9}{q_1}^3-\frac{8}{3}{q_1}^2q_2-\frac{8}{3}{q_1}{q_2}^2-\frac{8}{9}{q_2}^3\\
\dot{p_2} & = & -\frac{4}{3}q_1-\frac{8}{3}q_2-\frac{8}{9}\sqrt{3}q_1q_2-\frac{4}{9}\sqrt{3}{q_1}^2\\
          &-& \frac{16}{9}{q_2}^3-\frac{8}{3}{q_1}^2q_2-\frac{8}{3}{q_1}{q_2}^2-\frac{8}{9}{q_1}^3.         
\end{eqnarray*}
For a partial solution we get $p=p_1=p_2$, $q=q_1=q_2$ and we find $\psi (t)= -\displaystyle{\frac{3\sqrt{3}}{\sqrt{26}\sinh(2it)+1}}$.

 We get $\eta_1=dp_1$, $\eta_2=dp_2$, $\xi_1=dq_1$, $\xi_2=dq_2$, $\xi=\xi_1+\xi_2$ and we obtain for NVE
\begin{equation}
\label{eta}
\ddot{\xi}=-(4+\frac{8}{9}\sqrt{3}.\psi(t) +24.{\psi(t)}^2)\xi.
\end{equation}
We need to algebrize (\ref{eta}) with a standard change of variable $w=\sqrt{26}\sinh(2it)+1$. The result is
\begin{equation}
\label{EQ}
 \xi '' =r(w)\xi . 
 \end{equation}

Next we  use the Kovacic algorithm to show that (\ref{EQ}) has non Liouvillian solutions and the identity component of the Galois group of this equation is $SL(2,\, {\bf C})$. This means that the Hamiltonian system (\ref{L}) is non integrable with meromorphic  first integrals. This proves that in the case $H_2+H_3+H_4$ there is non additional meromorphic (holomorphic) first integral.
 The system $H_3$ is integrable, $H_2+H_3$ and $H_2+H_3+H_4$ are non integrable thats why we could conclude that the system $H_2+H_3+\dots + H_k $ is non integrable for each $k\ge 3$. This proves B).

{\bf Acknowledgements.} 
I acknowledge for the partial support by Sofia University Grant 80-10-215/2017.

{Georgi Georgiev \\
ggeorgiev3@fmi.uni-sofia.bg\\
 Faculty of Mathematics and Informatics, Sofia University, \\1164 Sofia, Bulgaria}
\end{document}